\documentstyle{article}

\newtheorem{theorem}{Theorem}

\newtheorem{definition}[theorem]{Definition}
\newtheorem{example}[theorem]{Example}
\newtheorem{remark}[theorem]{Remark}

\def\QED{\quad\blackslug\lower 8.5pt\null}

  \newcommand{\htimes}{%
    \mathbin{\mathsurround0pt \mathchoice
      {\hbox{\vrule\negthinspace$\times$}}%
      {\hbox{\vrule\negthinspace$\times$}}%
      {\hbox{\vrule\negthinspace$\scriptstyle\times$}}%
      {\hbox{\vrule\negthinspace$\scriptscriptstyle\times$}}%
    }}

\begin{document}

\noindent {\em Rendiconti del Seminario di Messina} \\
{\em Serie II, Vol. ? (1999)}

\vspace*{5mm}

\begin{center}
{\Large \bf  LIGHTLIKE  HYPERSURFACES ON A} 

\vspace*{2mm}

{\Large \bf FOUR-DIMENSIONAL MANIFOLD} 

\vspace*{2mm}

{\Large \bf  ENDOWED WITH A PSEUDOCONFORMAL}

\vspace*{2mm}

{\Large \bf  STRUCTURE OF SIGNATURE (2, 2)}

\vspace*{3mm}

{\large MAKS A. AKIVIS and  VLADISLAV V. GOLDBERG}

\vspace*{5mm}
{\em Dedicated to the memory of Professor Pasquale Calapso}

\end{center}

\vspace*{5mm}

{\footnotesize{\bf Abstract}. 
{\em The authors study the geometry of 
lightlike hypersurfaces on a four-dimensional 
manifold  $(M, c)$ endowed with 
a pseudoconformal structure $c = CO (2, 2)$. 
They  prove that a lightlike hypersurface 
$V \subset (M, c)$ 
bears a foliation formed by conformally invariant 
isotropic geodesics and two isotropic  distributions tangent to 
these geodesics, and that  these two distributions are  
integrable if and only if   $V$ is totally umbilical. 
The authors  also indicate how, using  singular points and 
singular submanifolds of a lightlike 
hypersurface $V \subset (M, c)$, 
to construct  an invariant normalization of $V$ intrinsically 
connected with $V$.} 
}
\vspace*{4mm}

\noindent
{\bf 1991 MS classification}: 53A30, 53B25.


\vspace*{5mm}

\setcounter{equation}{0}

\setcounter{section}{-1}

\section{INTRODUCTION}  

A four-dimensional pseudo-Riemannian manifold $(M, g)$  
with a metric quadratic form of signature $(3, 1)$ 
is a geometric model of the classic spacetime in 
general relativity. Its natural generalization is 
a  pseudo-Riemannian manifold $(M, g)$  
of dimension $n = \dim \; M$ 
with a nondegenerate quadratic form of arbitrary 
signature $(p, q), \; p + q = n.$ Such manifolds are considered 
in a construction of multidimensional models of spacetime and in 
the theory of superstrings.

Let $x$ be a point of  a manifold $(M, g)$, and $T_x (M)$  
the tangent space at the point $x$. For $p \cdot q > 0$, 
the quadratic form $g$ defines  a real {\em isotropic 
cone} $C_x$ in the space $T_x$.  
Its equation is $g (\xi, \xi) = 0, \; \xi \in T_x$. 
This cone is also called  the {\em light cone} 
or the {\em null cone}. 

A hypersurface $V \subset (M, g)$  is called {\em lightlike} 
if it is tangent to the cone 
$C_x$ at each  point $x \in V$. The lightlike 
hypersurfaces are also called {\em isotropic} or 
{\em null hypersurfaces}. On the manifold $(M, g)$,
 such hypersurfaces separate domains with different physical or 
geometric properties---they are models of physical or geometric 
horizons (see, for example,  [HE 73]). 

Many physical and geometric 
objects on a manifold $(M, g)$ are invariant under 
conformal transformations of the metric $g$, that is, 
under a passage from the metric $g$ to the metric $\widetilde{g} 
= \sigma g$, where $\sigma = \sigma (x)$ is a differentiable 
function such that $\sigma (x) \neq 0, \; x \in M$. 
 Examples of such objects are the light cones and 
the lightlike hypersurfaces. 
 Hence it is appropriate to study such objects  
not only on a pseudo-Riemannian manifold $(M, g)$ 
but also on a manifold $(M, c)$ endowed with a conformal 
structure $c = \{\sigma g\}$.

 Note that  lightlike hypersurfaces arose in 
the papers  of Duggal and Bejancu [DB 91] and [Be 96]  
(see also their book  [DB 96],  Section 4.7).
They considered them in a pseudo-Riemannian manifold of constant 
curvature $c$, and in particular, in pseudo-Euclidean spaces 
${\bf R}^4_1$ and ${\bf R}^4_2$. Lightlike hypersurfaces were 
also studied by Kupeli  [Ku 87] (see also his book [Ku 96], 
Section 4.4). He considered them in a (pseudo-)Riemannian space 
$(M, g)$ of constant sectional curvature.  
Lightlike hypersurfaces appeared in the paper Rosca [Ro 71] 
in which he studied a pair of lightlike hypersurfaces in 1-to-1 
correspondence in a Lorentz manifold. 

Note also that the totally umbilical lightlike hypersurfaces in 
Riemannian and pseudo-Riemannian manifolds $(M, g)$ were 
extensively studied by many authors. They considered 
their local and global properties. For example, in the papers 
[Y 75], [Ak 87], [Ra 87], and [Z 96] the authors found 
 necessary and sufficient conditions for a complete 
 spacelike hypersurface to be totally umbilical in $(M, g)$. 

The totally umbilical lightlike hypersurfaces in 
$(M, c)$ endowed with a conformal or pseudoconformal 
structure  were not yet studied extensively. 
In the papers [AG 99a] and [AG 99b] 
we have already studied the  geometry of 
 lightlike hypersurfaces $V$ 
on a manifold $(M, c)$ endowed with a conformal 
structure $c$ of Lorentzian signature $(n - 1, 1)$. 

In the present paper we 
consider  lightlike hypersurfaces on a 
four-dimensional manifold $(M, c)$ 
endowed with a conformal structure of ultrahyperbolic 
signature $(2, 2)$. We find some properties 
of the structure of such hypersurfaces and prove that 
they  bear two isotropic two-dimensional distributions
in addition to the fibration of isotropic geodesics. 
We also prove that  integrability of these distributions   
is necessary and sufficient 
for a lightlike  hypersurface to be totally umbilical. 
We constructed an invariant normalization of $V \subset (M, c)$ 
in a fourth differential neighborhood of a point of $V$. 
In the case in question, i.e., for $c = CO (2, 2)$, 
we were able not only to construct such normalization 
but also to find a foliation of canonical frames.


For our study of  lightlike hypersurfaces 
on a manifold $(M, c), \; \dim \; M = 4, \; \newline 
\mbox{{\rm sign}} \; c = (2, 2)$, we use the apparatus developed 
in [A 96] (see also [AG 96], Ch. 5). As far as we know, 
the lightlike hypersurfaces on such manifolds 
are  studied in the present paper for the first time.

\section{A MANIFOLD 
\protect\boldmath \( (M, c) \)\protect\unboldmath}

Consider a manifold $(M, c)$ endowed with a conformal 
structure $c$ of signature $(p, q), \; \dim \; M = n = p + q$. 
Let $x$ be an arbitrary point of $M, \; T_x (M)$ be its 
 the tangent space, and  $C_x \subset T_x (M)$ be the isotropic 
cone in $T_x (M)$. The space $T_x (M)$ can be compactified 
by adding the point at infinity $y$ and the isotropic 
cone $C_y$ with vertex $y$. After this enlargement, the space 
$T_x (M)$ becomes a  pseudoconformal 
space $(C^n_q)_x$ of the same signature $(p, q)$.  

 Under the Darboux mapping (see [AG 98] or [AG 96], Ch. 1), 
 the space $(C^n_q)_x$ will be 
mapped onto a hyperquadric $(Q_q^n)_x$ of a projective space  
 $P^{n+1}_x$ of dimension $n + 1$. In the space $P^{n+1}_x$, 
the hyperquadric $(Q_q^n)_x$ is  defined by the equation 
$(x, x) = 0$.

 We associate a family 
of projective local frames $\{A_0, A_i, A_{n+1}\}, \; 
i = 1, \ldots 
, n$, with this hyperquadric in such a way 
that $A_0 = x$ and $A_{n+1} = y$, where $x$ and $y$ are 
points of the hyperquadric $(Q_q^n)_x$ for which $(x, y) \neq 0$. 
This implies
$$
(A_0, A_0) = (A_{n+1}, A_{n+1}) = 0, \;\;(A_0, A_{n+1}) = -1
$$
The last condition is obtained by taking an appropriate 
normalization of the points $A_0$ and $A_{n+1}$. 
Here and in what follows the parentheses denote the scalar 
product 
with respect to the quadratic form occurring in the left-hand 
side 
of the equation of the hyperquadric $(Q^n_q)_x$.

Denote by 
$T_x$ and $T_y$ the  tangent hyperplanes to $(Q^n_q)_x$ in 
the points $x$ and $y$ and locate the points $A_i$ at the 
intersection of these hyperplanes, $A_i \in T_x \cap T_y, \; 
i = 1, \ldots, n$. Then we find that 
$$
(A_0, A_i)  = (A_{n+1}, A_i) = 0, \;\; (A_i, A_j) = g_{ij}, 
$$
where $\det (g_{ij}) \neq 0, \;\mbox{{\rm sign}}\; 
(g_{ij}) = (p, q)$ (see Figure 1).

\vspace*{50mm}
\begin{center}
Figure 1
\end{center}
Now 
\begin{equation}\label{eq:1}
\renewcommand{\arraystretch}{1.3}
\begin{array}{ll} 
(A_\xi, A_\eta) = (g_{\xi\eta}) = \pmatrix{0 & 0 & -1\cr
                         0 & g_{ij} & 0 \cr
						 -1 & 0 & 0\cr}, 
\;\;\;\; \xi, \eta = 0, 1, \ldots , n+1, 

\end{array}
\renewcommand{\arraystretch}{1}
\end{equation}
and the equation of the hyperquadric  $(Q^n_q)_x \subset 
P^{n+1}_x$ can be written in the form 
$$
(x, x) = g_{ij} x^i x^j - 2 x^0 x^{n+1} = 0.
$$
The family of projective frames we have constructed is 
called a {\em first order frame bundle associated with 
the manifold} $(M, c)$.

The isotropic cone $C_x$ is the intersection of 
the hyperquadric $(Q_q^n)_x$ and the tangent hyperplane 
$T_x$:
$$
C_x = T_x \cap (Q^n_q)_x,
$$ 
and it is defined by the equation
$$
g = g_{ij} \xi^i \xi^j = 0, \;\;\xi = (\xi^i) \in T_x.
$$
The group of transformations of the tangent space $T_x (M)$ 
preserving the invariant cone 
$C_x$ is the group $G = {\bf SO} (p, q) \times {\bf H}$,
where ${\bf SO} (p, q)$ is a special pseudoorthonormal 
group of signature $(p, q)$ and ${\bf H} = {\bf R}^+$ 
is the group of homotheties. 

It follows from relations (1) that 
the equations of infinitesimal displacement  
in the first-order  frame bundle have the form
\begin{equation}\label{eq:2}
\renewcommand{\arraystretch}{1.3}
\left\{
\begin{array}{lll}
dA_0 = \omega_0^0  A_0 + \!\!\!\!& \!\!\!  \omega_0^i  A_i, &\\
dA_i = \omega_i^0  A_0 + \!\!\!\! & \!\!\! \omega_i^j  A_j   
\!\!\!\!& + \omega_i^{n+1} A_{n+1}, \\
dA_{n+1} = &\!\!\! \omega_{n+1}^i  A_i  \!\!\!\! &
- \omega^0_0 A_{n+1}, 
\end{array} 
\right.
\renewcommand{\arraystretch}{1}
\end{equation}
where 
\begin{equation}\label{eq:3}
\renewcommand{\arraystretch}{1.3}
\begin{array}{ll} 
\omega_0^i = \omega^i \\
 \omega_i^{n+1} = g_{ij}  \omega^j, \;\; 
 \omega^i_{n+1} = g^{ij}  \omega_j^0, 
\\
d g_{ij} - g_{ik} \omega_j^k  -  g_{kj} \omega_i^k  = 0,
\renewcommand{\arraystretch}{1}
\end{array}
\end{equation}
and $g^{ij}$ is the inverse tensor of the tensor $g_{ij}$, 
i.e., $g^{ik} g_{kj} = \delta^i_j$. Note that the tensor $g_{ij}$ and the 1-forms $\omega^i$ 
are defined in a first-order neighborhood of a point 
$x \in (M, c)$,   the 1-forms $\omega^0_0$ 
and $\omega_j^i$ are defined in its a second-order 
neighborhood, and the 1-forms $\omega^0_i$ are defined in its  
third-order neighborhood.

The 1-forms $\omega^i$ define a displacement of a point $A_0$ 
and consequently of a frame $\{A_0, A_i, A_{n+1}\}$ 
along the manifold $M$. This is the reason that they are 
called the {\em basis forms}. For $\omega^i = 0$, 
equations (2) take the form 
\begin{equation}\label{eq:4}
\renewcommand{\arraystretch}{1.3}
\left\{
\begin{array}{lll}
\delta A_0 = \pi_0^0  A_0, &  &\\
\delta A_i = \pi_i^0  A_0 + & \!\!\!\! \pi_i^j  A_j,& \\
\delta A_{n+1} = &\!\!\! \pi_{n+1}^i  A_i 
                          - &\!\!\!\! \pi^0_0 A_{n+1}. 
\end{array} 
\right.
\renewcommand{\arraystretch}{1}
\end{equation}
Here $\delta$ is the symbol of differentiation for 
$\omega^i = 0$, i.e., with respect to 
the fiber parameters of the frame bundle, and 
$\pi_\eta^\xi = \omega_\eta^\xi (\delta)$. Formulas 
(4) define admissible transformations in a fiber of 
a  first-order  frame bundle. These transformations form the 
group  $G' = G \htimes {\bf T} (n)$ that is obtained by  
a differential prolongation of the group $G$ acting 
in the space $T_x (M)$. Here ${\bf T} (n)$ is a subgroup 
of the group $G'$ which is isomorphic to the group 
of parallel translations, and the symbol $\htimes$ 
is the semidirect product (see [AG 96], Ch. 4). Equations 
(4) show that the group $G'$  is isomorphic to the subgroup 
of the group of pseudoconformal transformations of the space $C_q^n$ keeping invariant the point $x = A_0$ of this space.

\section{THE TENSOR OF CONFORMAL \\ CURVATURE OF   
A MANIFOLD \\ 
\protect\boldmath \( (M, c), \; c = CO (2, 2) \) 
\protect\unboldmath}

As was proved in Ch. 4 of the book [AG 96], 
the structure equations of the manifold 
$(M, c)$ endowed with a conformal structure 
of an arbitrary signature $(p, q)$ have the form
\begin{equation}\label{eq:5}
\renewcommand{\arraystretch}{1.3}
\left\{
\begin{array}{ll}
d \omega^i = \omega_0^0 \wedge \omega^i 
+ \omega^j \wedge \omega^i_j,\\ 
 d \omega_0^0 = \omega^i \wedge \omega^0_i\\ 
d \omega^i_j = \omega_j^0 \wedge \omega^i + 
 \omega^k_j \wedge \omega_k^i +  \omega^{n+1}_j \wedge 
 \omega_{n+1}^i + C^i_{jkl} \omega^k \wedge \omega^l,\\
d \omega_i^0 = \omega^0_i \wedge \omega_0^0
 + \omega_i^j \wedge \omega^0_j 
+ C_{ijk}  \omega^j  \wedge  \omega^k.
\end{array}
\right.
\renewcommand{\arraystretch}{1}
\end{equation}
Here the quantities $C^i_{jkl}$ are defined in a 
third-order neighborhood of a point $x \in (M, c)$ 
and form the {\em tensor of conformal curvature}, also 
called the {\em Weyl tensor}. Denote it by the letter 
$C$, where $C = (C_{jkl}^i)$.

The quantities $C_{ijk}$ are defined in a 
fourth-order neighborhood of a point $x \in (M, c)$. 
and for $n \geq 4$, they do not form a tensor. 
Denote the object  $C_{ijk}$ by $C'$, i.e. 
$C' = (C_{ijk})$. The reason for this notation is that 
 for $n \geq 4$, this object is expressed 
in terms of the covariant derivatives of the tensor $C$. 
For $n \geq 4$, the condition $C = 0$ implies $C' = 0$, and 
a manifold $(M, c)$ is conformally flat, i.e., it is locally 
equivalent to a conformal space $C^n_q$.

For $n = 3$, the tensor $C = (C^i_{jkl})$ is identically 
equal to 0, and the curvature of the space is defined by the 
object $C' = (C_{ijk})$ which in this case becomes a tensor. 
In what follows, we will assume that $n \geq 4$.

The components of the tensor $C$ and the object $C'$ satisfy the 
equations
\begin{equation}\label{eq:6}
\renewcommand{\arraystretch}{1.3}
\left\{
\begin{array}{ll}
C_{ijkl} = g_{im}  C^m_{jkl},\\
C_{ijkl} = - C_{jikl} = - C_{ijlk},  & C_{ijkl} = C_{klij} ,\\ 
C_{ijkl} + C_{iklj} + C_{iljk} = 0,\\
C^i_{jki} = 0, \;\;C_{ijk} = - C_{ikj}.
\end{array} 
\right.
\renewcommand{\arraystretch}{1}
\end{equation}

Note that we will use the notation $C$ not only for the 
tensor $C^i_{jkl}$ but also for the tensor $C_{ijkl}$.

\section{ISOTROPIC FRAMES FOR   \\
\protect\boldmath \( (M, c), \; c = CO (2, 2) \) 
\protect\unboldmath}  

Consider four-dimensional manifold $(M, c)$ endowed 
with a pseudoconformal structure $c = CO (2, 2)$. The 
fundamental quadratic form of this signature 
is reduced to the form
\begin{equation}\label{eq:7}
g = 2 (\xi^2 \xi^3 - \xi^1 \xi^4).
\end{equation}
It follows that the matrix of its coefficients is 
\begin{equation}\label{eq:8}
(A_i, A_j) = (g_{ij}) = \pmatrix{0 & 0 & 0& -1\cr
                         0 & 0 & 1 & 0 \cr
                         0 & 1 & 0 & 0 \cr
					-1 & 0 & 0 & 0\cr}.
\end{equation}

Note that we changed the signs of components 
 of the tensor $g_{ij}$ in comparison with the book [AG 96] 
and the paper [A 96].

The isotropic cone $C_x \subset T_x (M)$ is defined by the 
equation
$$
\xi^2 \xi^3 - \xi^1 \xi^4 = 0.
$$
We will clarify a structure of this cone. The last equation  
can be written in two different ways:
$$
\frac{\xi^1}{\xi^3} = \frac{\xi^2}{\xi^4} = - \lambda, \;\;\;\;
\frac{\xi^1}{\xi^2} = \frac{\xi^3}{\xi^4} = - \mu.
$$ 
It follows that the cone $C_x$ carries two families of 
two-dimensional plane generators. These are defined by the equations 
\begin{equation}\label{eq:9}
 \xi^1 + \lambda  \xi^3 = 0, \;\;  \xi^2 + \lambda  \xi^4 = 0
\end{equation} 
and
\begin{equation}\label{eq:10}
 \xi^1 + \mu  \xi^2 = 0, \;\;  \xi^3 + \mu  \xi^4 = 0.
\end{equation}
The 2-planes defined by equations (9) are called 
$\alpha$-generators, and the 2-planes defined by 
equations (10) are called  $\beta$-generators 
of the cone $C_x$. 

\vspace*{1.8in}
\begin{center}
Figure 2
\end{center}

Under the projectivization of the tangent space 
$T_x (M)$ with center at a point $x = A_0$, 
there corresponds a ruled quadric $P C_x$ 
 for the cone $C_x$ where $PC_x$ belongs to  a three-dimensional 
projective space $P^3_x = PT_x (M)$. With respect to the 
frame $\{\widetilde{A}_i\}$, where $\widetilde{A}_i = P A_i$, 
the quadric $PC_x$ is defined by the same equation (7) 
(see Figure 2).

For the conformal structure $CO (2, 2)$, 
the group of transformations of the tangent space $T_x (M)$ 
preserving the invariant cone is split into three subgroups:
 $G = {\bf SL} (2)  \times {\bf SL} (2) \times {\bf H}$. 
The first two of these subgroups transfer the families 
of $\alpha$- and $\beta$-generators of the cones $C_x$ 
into themselves and are isomorphic to the group 
of projective transformations on a projective straight line 
$P^1$. As in the general case, the third subgroup is the group of 
homotheties. 

On the manifold $(M, c)$, the isotropic $\alpha$- and 
$\beta$-generators of the cone $C_x$ 
form two fiber bundles $E_\alpha$ and $E_\beta$ with 
the common base $M$. The fibers of $E_\alpha$ and $E_\beta$ 
are the families of $\alpha$- and $\beta$-generators of the cones 
$C_x$. By (9) and (10), these fibers are parameterized by means of 
nonhomogeneous projective parameters $\lambda$ and $\mu$ and 
are isomorphic to real projective straight lines 
${\bf R} P_\alpha$ and ${\bf R} P_\beta$. Thus  the 
fiber bundles $E_\alpha$ and $E_\beta$ can be written as 
$E_\alpha = (M, {\bf R} P_\alpha)$ and 
$E_\beta = (M, {\bf R} P_\beta)$. These fiber bundles are {\em 
real twistor fibrations} similar to those introduced 
on four-dimensional manifolds of Lorentzian signature $(3, 1)$ by 
Penrose (see, for example, [PR 86]). 

Consider $\alpha$- and $\beta$-generators of the cone $C_x$. 
For $\omega^i = 0$ (i.e., for fixed principal parameters), they 
are defined by equations (9) and (10). One can easily prove that 
these generators intersect one another in an isotropic straight 
line connecting the point $A_0 = x$ with the point 
\begin{equation}\label{eq:11}
B = \lambda \mu A_1 - \lambda A_2 - \mu A_3 + A_4,
\end{equation} 
and that they belong to a three-dimensional subspace of the 
space $T_x (M)$ defined by the equation
\begin{equation}\label{eq:12}
\xi^1 + \mu \xi^2 + \lambda \xi^3 + \lambda \mu \xi^4 = 0.
\end{equation} 
This subspace is tangent to the isotropic cone $C_x$ along its 
generator $A_0 B$ and is also called isotropic.

In the space $T_x (M)$, we specialize our moving frame 
in such a way that its vertex $A_1$ coincides with the point $B$ and 
the isotropic straight line $A_0 B$ coincides with the straight 
line $A_0 A_1$. Then the nonhomogeneous projective parameters 
$\lambda$ and $\mu$ occurring in equations (9) and (10) 
for isotropic 2-planes $A_0 A_1 A_2$ and $A_0 A_1 A_3$ 
become $\infty \;$, $\lambda = \infty, \; \mu = \infty$, 
and equations (9) and (10) take the form
\begin{equation}\label{eq:13}
 \xi^3 = 0, \;\;   \xi^4 = 0
\end{equation} 
and
\begin{equation}\label{eq:14}
  \xi^2 = 0, \;\;   \xi^4 = 0.
\end{equation}
The equation of the isotropic subspace (12) containing these 
isotropic $\alpha$- and $\beta$-generators becomes 
\begin{equation}\label{eq:15}
   \xi^4 = 0.
\end{equation}

\section{THE STRUCTURE EQUATIONS OF \\ A MANIFOLD 
\protect\boldmath \( (M, c) \) \protect\unboldmath}   

For a manifold $(M, c)$ endowed with a conformal 
structure $c = CO (2, 2)$, in the isotropic frame bundle, 
equations (8) and the last equation of (3) imply that 
\begin{equation}\label{eq:16}
\renewcommand{\arraystretch}{1.3}
\left\{
\begin{array}{ll}
\omega_1^4 = \omega_2^3 = \omega_3^2 = \omega_4^1 = 0,\\
\omega_2^4 = \omega_1^3, \;\;\;\;\;\;\;\; 
\omega_4^2 = \omega_3^1, \\
\omega_3^4 = \omega_1^2, \;\; \;\;\;\;\;\; 
\omega_4^3 = \omega_2^1, \\
\omega_1^1 + \omega_4^4 = 0, \;\;  \omega_2^2 + \omega_3^3 = 0.
\end{array} 
\right.
\renewcommand{\arraystretch}{1}
\end{equation}
Thus on the manifold $(M, c)$, among 
the forms $\omega_i^j$  only  the forms  $\omega_1^2, \; 
\omega_2^1, \; \omega_1^3, \; \omega_3^1, \; \omega_1^1,$ and 
$\omega_2^2$ are independent. Hence 
on such a manifold $(M, c)$,  the structure equations (5) 
take the form
\begin{equation}\label{eq:17}
\left\{
\renewcommand{\arraystretch}{1.3}
\begin{array}{ll}
d\omega^1 = (\omega_0^0 - \omega_1^1) \wedge \omega^1 + \omega^2 
\wedge \omega_2^1 + \omega^3 \wedge \omega_3^1, \\
d\omega^2 = (\omega_0^0 - \omega_2^2) \wedge \omega^2 + \omega^1 
\wedge \omega_1^2 + \omega^4 \wedge \omega_3^1, \\
d\omega^3 = (\omega_0^0 + \omega_2^2) \wedge \omega^3 + \omega^1 
\wedge \omega_1^3 + \omega^4 \wedge \omega_2^1, \\
d\omega^4 = (\omega_0^0 + \omega_1^1) \wedge \omega^4 + \omega^2 
\wedge \omega_1^3 + \omega^3 \wedge \omega_1^2,
\end{array}
\right.
\renewcommand{\arraystretch}{1}
\end{equation}
\begin{equation}\label{eq:18}
 d \omega_0^0 = \omega^1 \wedge \omega_1^0 
+ \omega^2 \wedge \omega_2^0 
+ \omega^3 \wedge \omega_3^0 + \omega^4 \wedge \omega_4^0,
\end{equation}
\begin{equation}\label{eq:19}
\renewcommand{\arraystretch}{1.3}
\left\{
\begin{array}{rl}
d \omega^3_1 = & \omega_1^0 \wedge \omega^3 
+ \omega_2^0 \wedge \omega^4  
+  (\omega_1^1 + \omega_2^2)  \wedge   \omega^3_1 + \Omega_1^3.\\
d (\omega^1_1 + \omega^2_2) = &
 \omega_1^0 \wedge \omega^1 + \omega_2^0 \wedge \omega^2  
- \omega_3^0 \wedge \omega^3 
- \omega_4^0 \wedge \omega^4 \\&
+ 2 \omega^3_1 \wedge \omega_3^1  + \Omega_1^1 + \Omega_2^2, \\
d \omega^1_3= &\omega^0_3 \wedge \omega^1 
+ \omega^0_4 \wedge \omega^2  
+    \omega^1_3  \wedge (\omega_1^1 + \omega_2^2)  
+ \Omega_3^1,
\end{array}
\right.
\renewcommand{\arraystretch}{1}
\end{equation} 
and 
\begin{equation}\label{eq:20}
\left\{
\renewcommand{\arraystretch}{1.3}
\begin{array}{rl}
d \omega^2_1 =& \omega^0_1 \wedge \omega^2 
+ \omega^0_3 \wedge \omega^4  + (\omega^1_1 - \omega_2^2) 
 \wedge \omega_1^2 + \Omega_1^2, \\
d (\omega^1_1 - \omega^2_2)  = & 
 \omega^0_1 \wedge \omega^1 - \omega^0_2 \wedge \omega^2  
+ \omega^0_3 \wedge \omega^3 - \omega^0_4 \wedge \omega^4 
\\&+ 2 \omega^2_1 \wedge \omega_2^1 
+ \Omega_1^1 - \Omega_2^2,\\
d \omega^1_2 =& \omega_2^0 \wedge \omega^1 
+ \omega^0_4 \wedge \omega^3  + \omega^1_2  
 \wedge (\omega_1^1 - \omega_2^2) + \Omega_2^1. 
\end{array}
\right.
\renewcommand{\arraystretch}{1}
\end{equation} 

Equations (19) prove that the 1-forms $\omega_1^3, \;
\omega_1^1 + \omega_2^2,$ and $\omega_3^1$ are fiber 
forms on the isotropic fiber bundle 
$E_\alpha$, and the 2-forms $\Omega_1^3, \;
\Omega_1^1 + \Omega_2^2,$ and $\Omega_3^1$ are 
the curvature forms of this  fiber bundle. 
Similarly, it follows from equations (20) that 
the 1-forms $\omega_1^2, \;
\omega_1^1 - \omega_2^2,$ and $\omega_2^1$ are fiber 
forms on the isotropic fiber bundle 
$E_\beta$, and the 2-forms $\Omega_1^2, \;
\Omega_1^1 - \Omega_2^2,$ and $\Omega_2^1$ are 
the curvature forms of $E_\beta$. 

Since each of the indices $i, j, k,$ and $l$ takes only four  
values 1, 2, 3, 4, equations (6) imply that 
the tensor of conformal curvature $C_{ijkl}$ has 21 essential nonvanishing components that 
satisfy 11 independent conditions: 
\begin{equation}\label{eq:21}
\left\{
\renewcommand{\arraystretch}{1.3}
\begin{array}{ll}
C_{1234} - C_{1324} + C_{1423} = 0, \\
C_{1224} = C_{1334} = C_{1213} = C_{2434} = 0. \\
C_{1314} - C_{1323} = C_{1424} - C_{2324} = 0, \\
C_{1214} + C_{1223} = C_{1434} + C_{2334} = 0, \\
C_{1414} = C_{2323} = C_{1234} + C_{1324}.
\end{array}
\renewcommand{\arraystretch}{1}
\right.
\end{equation}
Hence the  tensor $C_{ijkl}$  has  10 
independent components in all:
\begin{equation}\label{eq:22}
\left\{
\renewcommand{\arraystretch}{1.3}
\begin{array}{ll}
C_{1212} = a_0, \;\;  C_{1214} = a_1, \;\; C_{1234} = a_2, \;\;
C_{1434} = a_3, \;\; C_{3434} = a_4, \\
C_{1313} = b_0, \;\; C_{1314} = b_1, \;\;C_{1324} = b_2, \;\;
 C_{1424} = b_3, \;\; C_{2424} = b_4.
\end{array}
\renewcommand{\arraystretch}{1}
\right.
\end{equation}
As a result, the curvature forms of the 
fiber bundles  $E_\alpha$ and $E_\beta$ can be written as 
\begin{equation}\label{eq:23}
\renewcommand{\arraystretch}{1.3}
\begin{array}{rlrl}
\Omega_1^3& \!\!\!\!= & 
-2 &\!\!\!\! [a_0 \omega^1 \wedge \omega^2 
+ a_1 (\omega^1 \wedge \omega^4 - \omega^2 \wedge \omega^3) 
+ a_2 \omega^3 \wedge \omega^4], \\
\Omega_1^1 + \Omega_2^2& \!\!\!\!
= & \!\!\!\! + 4&\!\!\!\! [a_1 \omega^1 \wedge \omega^2 
+ a_2 (\omega^1 \wedge \omega^4 - \omega^2 \wedge \omega^3) 
+ a_3 \omega^3 \wedge \omega^4], \\
\Omega_3^1&\!\!\!\! = & \!\!\!\! + 
2&\!\!\!\! [a_2 \omega^1 \wedge \omega^2 
+ a_3 (\omega^1 \wedge \omega^4 - \omega^2 \wedge \omega^3) 
+ a_4 \omega^3 \wedge \omega^4], 
\end{array}
\renewcommand{\arraystretch}{1}
\end{equation}
and   
\begin{equation}\label{eq:24}
\renewcommand{\arraystretch}{1.3}
\begin{array}{rlrl}
\Omega_1^2& \!\!\!\!= &-2& \!\!\!\! [b_0 \omega^1 \wedge \omega^3 
+ b_1 (\omega^1 \wedge \omega^4 + \omega^2 \wedge \omega^3) 
+ b_2 \omega^2 \wedge \omega^4], \\
\Omega_1^1 - \Omega_2^2&\!\!\!\! 
= & \!\!\!\! + 4& \!\!\!\![b_1 \omega^1 \wedge \omega^2 
+ b_2 (\omega^1 \wedge \omega^4 + \omega^2 \wedge \omega^3) 
+ b_3 \omega^2 \wedge \omega^4], \\
\Omega_2^1&\!\!\!\! =&\!\!\!\! 
+ 2& \!\!\!\![b_2 \omega^1 \wedge \omega^3 
+ b_3 (\omega^1 \wedge \omega^4 + \omega^2 \wedge \omega^3) 
+ b_4 \omega^2 \wedge \omega^4]. 
\end{array}
\renewcommand{\arraystretch}{1}
\end{equation}

From equations (23) and (24) it follows 
the  tensor $C = (C_{ijkl})$ of conformal curvature 
is split into two subtensors $A$ and $B, \; C = A + B$, where 
$$
A = \{a_u\}, \;\; B = \{b_u\}, \;\;\;\; u = 0, 1, 2, 3, 4.
$$
These  are the curvature tensors of 
 the fiber bundles  $E_\alpha$ and $E_\beta$.

If one of subtensors $A$ or $B$ vanishes, then 
a manifold $(M, c)$ is called {\em conformally semiflat}. 
In this case  the 
fiber bundle  $E_\alpha$ (respectively, $E_\beta$) admits 
a three-parameter family of two-dimensional integral surfaces 
$V_\alpha$ (respectively, $V_\beta$). 
 
If both subtensors $A$ and $B$ vanish, then 
the tensor $C$ also vanishes. 
In this case a manifold $(M, c)$ becomes conformally flat 
and is locally equivalent to a pseudoconformal space $C^4_2$. 
Under the Darboux mapping, 
a hyperquadric $Q^4_2$ of a projective space $P^5$ corresponds 
for the space $C^4_2$.  
Under this mapping, 
two-dimensional plane generators of the hyperquadric 
$Q^4_2$.  correspond
for two-dimensional integral surfaces $V_\alpha$ and $V_\beta$ of the 
fiber bundles  $E_\alpha$ and $E_\beta$.  

\section{PRINCIPAL ISOTROPIC BIVECTORS} 
Suppose that $\xi$ and $\eta$ are vectors of the space $T_x (M)$, 
and $p = \xi \wedge \eta$ is a bivector defined by $\xi$ and 
$\eta$. The coordinates of $p$ are 
$$
p^{ij} = \xi^{[i} \eta^{j]} = \frac{1}{2} (\xi^i \eta^j - \xi^j \eta^i), \;\; p^{ij} = - p^{ji}.
$$
The tensor of conformal curvature  $C = (C_{ijkl})$ allows us to 
define the {\em relative conformal curvature of the bivector $p$}:
\begin{equation}\label{eq:25}
C (p) = C_{ijkl} p^{ij} p^{kl}.
\end{equation}

Let us find relative conformal curvatures of the bivectors 
$p_\lambda$ and $p_\mu$ defined by 
$\alpha$- and $\beta$-generators of the isotropic cone $C_x$ of 
the manifold $(M, c)$. From equations (9) it follows that 
the vectors $\xi_\lambda$ and $\eta_\lambda$ defining 
the bivector $p_\lambda$ are defined by the formulas 
$$
\xi_\lambda = e_3 - \lambda e_1, \;\; \eta_\lambda = e_4 
- \lambda e_2.
$$
As a result, the coordinates of the bivector $p_\lambda$ are
$$
p^{12} = \lambda^2, \;\; p^{13} = 0, \;\; 
p^{14} = - \lambda, \;\; p^{23} = \lambda, \;\; 
p^{34} = 1, \;\; p^{42} = 0.
$$
Substituting these values of coordinates $p^{ij}$ 
into formula (25) 
and applying relations (22), we find that 
$$
 \frac{1}{4} C (p_\lambda) = 
a_0 \lambda^4 - 4 a_1 \lambda^3 + 6 a_2 \lambda^2 
- 4 a_3 \lambda + a_4.
$$
Since the right-hand side of the last equation contains only 
the components of the curvature tensor $A$ of the 
isotropic fiber bundle $E_\alpha$, this formula can be written as 
\begin{equation}\label{eq:26}
 \frac{1}{4} A (p_\lambda) = 
a_0 \lambda^4 - 4 a_1 \lambda^3 + 6 a_2 \lambda^2 
- 4 a_3 \lambda + a_4.
\end{equation}

Similarly, the  bivector $p_\mu$ is defined by the vectors 
$$
\xi_\mu = e_2 - \mu e_1, \;\; \xi_\mu = e_4 - \mu e_3,
$$
and its coordinates  are
$$
p^{12} = 0, \;\; p^{13} = \mu^2, \;\; 
p^{14} = - \mu, \;\; p^{23} = -\mu, \;\; 
p^{34} = 0, \;\; p^{42} = -1.
$$
This implies that 
the relative conformal curvatures of the bivector  $p_\mu$ is
\begin{equation}\label{eq:27}
\frac{1}{4} B (p_\mu) = 
b_0 \mu^4 - 4 b_1 \mu^3 + 6 b_2 \mu^2 - 4 b_3 \mu + b_4.
\end{equation}

The isotropic bivectors whose relative conformal curvature
vanishes are called the {\em principal isotropic bivectors}. 
By (26) and (27), the values of parameters $\lambda$ and $\mu$ 
defining such bivectors satisfy the algebraic equations
\begin{equation}\label{eq:28}
a_0 \lambda^4 - 4 a_1 \lambda^3 + 6 a_2 \lambda^2 
- 4 a_3 \lambda + a_4 = 0
\end{equation}
and 
\begin{equation}\label{eq:29}
b_0 \mu^4 - 4 b_1 \mu^3 + 6 b_2 \mu^2 - 4 b_3 \mu + b_4 = 0.
\end{equation}
Thus in the general case 
 the isotropic cone $C_x$ bears four principal 
$\alpha$-generators and the same number of principal 
$\beta$-generators if we count each of these generators as many times as its multiplicity.

On a manifold $(M, c)$, the principal isotropic 
bivectors form four 
principal $\alpha$-distributions and the same number of 
principal $\beta$-distributions. It was proved in [A 96] 
that {\em if $\lambda$ is a multiple root of equation 
$(28)$, then the principal $\alpha$-distribution defined by 
this root 
is integrable.} In the same way {\em if $\mu$ is a multiple 
root of equation $(29)$, then the principal 
$\beta$-distribution defined 
by this root is integrable.}

Suppose that $\lambda$ and $\mu$ are two fixed roots of equations 
(28) and (29), respectively, and $p_\lambda$ and $p_\mu$ are 
the principal isotropic distributions defined by these two 
roots. By means of a frame transformation indicated at the end
 of Section {\bf 3}, the values of parameters $\lambda$ and $\mu$
 can be made to equal $\infty$, $\; \lambda = \infty, \;\mu 
= \infty$. As a result, by (9) and (10), we find that these two 
distributions are defined by the following two systems of 
equations:
\begin{equation}\label{eq:30}
\omega^3 = 0, \;\; \omega^4 = 0
\end{equation}
and
\begin{equation}\label{eq:31}
\omega^2 = 0, \;\; \omega^4 = 0.
\end{equation}
Moreover, equations (28) and (29) become 
\begin{equation}\label{eq:32}
- 4 a_1 \lambda^3 + 6 a_2 \lambda^2 
- 4 a_3 \lambda + a_4 = 0
\end{equation}
and 
\begin{equation}\label{eq:33}
 - 4 b_1 \mu^3 + 6 b_2 \mu^2 - 4 b_3 \mu + b_4 = 0.
\end{equation}

If the coefficient $a_1$ in equation (32) vanishes, then 
the root $\lambda = \infty$ of this equation is a multiple 
root, and as a result, the principal distribution (30) 
defined by this root is integrable. Two-dimensional 
integral surfaces $V_\alpha$ of this distribution form 
an isotropic fiber bundle on the manifold $(M, c)$. 
Similarly, if the coefficient $b_1$ in equation (33) vanishes, 
then the root $\mu = \infty$ of this equation is a multiple 
root and the principal distribution (31) 
defined by this root is integrable. In addition, the two-dimensional 
integral surfaces $V_\beta$ of this distribution form 
an isotropic fiber bundle on the manifold $(M, c)$.

\section{LIGHTLIKE HYPERSURFACES ON \\
  \protect\boldmath \( (M, c), 
\; c = CO (2, 2) \) \protect\unboldmath}

As we already said in the introduction, a hypersurface 
$V$ on a manifold $(M, c), \newline \dim V = 3$, is said to be 
{\em lightlike} if its tangent subspace $T_x (V)$ is 
{\em tangent to the isotropic cone $C_x$}, i.e., 
this subspace is isotropic. The aim of this paper is 
to study the geometry of lightlike  hypersurfaces 
on a manifold $(M, c)$, where $c = CO (2, 2)$.

With a point $x$ of a lightlike hypersurface $V$, we 
associate a moving frame in such a way 
that its vertex $A_0$ coincide with $x \in V, \;  A_0 = x$, 
the points $A_1, A_2$, and $A_3$ belong to the tangent subspace 
$T_x (V)$, and the point $A_1$ belongs to the  
isotropic straight line along which the subspace $T_x (V)$ is 
tangent to the isotropic cone $C_x$. The subspace $T_x (V)$ 
contains two isotropic $\alpha$- and $\beta$-planes 
intersecting one another along the straight line $A_0 A_1$. 
Thus the 2-plane $A_0 \wedge A_1 \wedge A_2$ is an 
$\alpha$-generator of the cone $C_x$, and the 2-plane $A_0 \wedge A_1 \wedge A_3$ is its $\beta$-generator. 

We place points $A_2$ and $A_3$ of our moving frame to 
these two planes and normalize them by the condition 
$(A_2, A_3) = 1$. The subspace $A_0 \wedge  A_2 \wedge  A_3$ 
 is called the {\em screen subspace} and is denoted by $S_x\;,$
$S_x = A_0 \wedge  A_2 \wedge  A_3 \subset T_x (V)$. 
Further we take a point $A_4$ on the isotropic cone $C_x$ 
in such a way that the subspace  $A_0 \wedge  A_1 \wedge  A_4$ 
is conjugate to the subspace $S_x$ with respect to 
the cone $C_x$. In addition, we normalize the points $A_1$ 
and  $A_4$ by the condition $(A_1, A_4) = - 1$. 

A straight line $N_x = A_0 \wedge A_4$ 
does not belong to the tangent subspace 
$T_x (V)$. This line  is called a {\em normalizing straight line}. 
Its location is uniquely determined by the location of the subspace $S_x$.

The matrix of scalar products of the points $A_i. \; 
i = 1, 2, 3, 4,$ now has the form (8). 

the family of frames we have constructed is called 
a {\em family of first-order frames associated with 
a point $x$ of a lightlike hypersurface} $V \subset (M, c)$.

We will now find the equations of a 
bundle of first-order isotropic frames  associated with a 
lightlike hypersurface  $V$. 
Since its tangent subspace $T_x (V) = A_0 \wedge 
A_1 \wedge A_2 \wedge A_3$, with respect to this frame bundle 
the equation of $V$ is 
\begin{equation}\label{eq:34}
\omega^4 = 0,
\end{equation}
and as a result, we have
\begin{equation}\label{eq:35}
dA_0 = \omega^0_0 A_0 + \omega^1 A_1 + \omega^2 A_2 
+ \omega^3 A_3. 
\end{equation}
The 1-forms 
$\omega^1, \omega^2$, and $\omega^3$ are independent. They 
are {\em basis forms} of the frame bundle in question and of 
the hypersurface $V$.

Equations 
\begin{equation}\label{eq:36}
\omega^2 = \omega^3 = 0 
\end{equation}
define on $V$ a foliation 
formed by isotropic lines. As was proved in [AG 99b], 
these lines are isotropic geodesics for all pseudo-Riemannian 
metrics $g$ compatible with the conformal structure 
$c = CO (2, 2)$ on the manifold $(M, c)$. 

We will assume that the isotropic geodesics defined 
by equations (36) can be prolonged indefinitely on a hypersurface 
$V$. In this case each of these geodesics bears the geometry 
of a projective straight line $P^1$, and a hypersurface $V$ is the image of the product $M^2 \times P^1$ under its 
differentiable mapping $f$ into the manifold $(M, c)$: 
$V = f (M^2 \times P^1), \; f: M^2 \times P^1 \rightarrow 
(M, c)$. 

Equation $\omega^3 = 0$ defines on $V$  a fibration of isotropic 
$\alpha$-planes $A_0 \wedge 
A_1 \wedge A_2$, and equation $\omega^2 = 0$ defines on $V$  a 
fibration of isotropic $\beta$-planes $A_0 \wedge 
A_1 \wedge A_3$ (cf. these two equations 
with equations (30) and (31)).

In an isotropic frame bundle, the first fundamental form $I$ of 
a lightlike hypersurface $V \subset (M, c)$ becomes
\begin{equation}\label{eq:37}
I = g|_V = (dA_0, dA_0) = 2 \omega^2 \omega^3.
\end{equation}
This form is of rank  2 and of signature $(1, 1)$, and its 
coefficients form the matrix
\begin{equation}\label{eq:38}
(g_{ab}) = \pmatrix{0 & 1 \cr 
                    1  & 0 \cr}, \; a, b = 2, 3.
\end{equation}

\section{SINGULAR POINTS AND TOTALLY  \\ UMBILICAL 
HYPERSURFACES}

By the last equation of (17), 
exterior differentiation of equation (34) (the basic equation 
of a lightlike hypersurface $V$) leads to the following exterior 
quadratic equation:
\begin{equation}\label{eq:39}
\omega^2 \wedge \omega_1^3 + \omega^3 \wedge \omega_1^2 = 0.
\end{equation}
Applying Cartan's lemma to this equation, we that  
\begin{equation}\label{eq:40}
\renewcommand{\arraystretch}{1.3}
\left\{
\begin{array}{ll}
\omega_1^3 = \lambda_{22} \omega^2 +  \lambda_{23}   
\omega^3, \\
\omega_1^2 = \lambda_{32} \omega^2 +  \lambda_{33} 
\omega^3,
\end{array} 
\right.
\renewcommand{\arraystretch}{1}
\end{equation}
where $\lambda_{23} = \lambda_{32}$.

By means 
of the Cartan test (see [BCGGG 91] and cf. [AG 99a]), one can prove 
that lightlike hypersurfaces $V \subset (M, c)$, where 
$c = CO (2, 2)$, exist and depend on a function of two 
variables.

Differentiating equation (35), we obtain
\begin{equation}\label{eq:41}
d^2 A_0 \equiv (\omega^2 \omega_2^4 + \omega^3 \omega_3^4) A_4 
+  (\omega^2 \omega_2^5 + \omega^3 \omega_3^5) A_5 
     \pmod{T_x (V)}.
\end{equation}
But by (3) and (8) we have
$$
\omega_2^5 = \omega^3, \;\;  \omega_3^5 = \omega^2, \;\;
\omega_2^4 = \omega_1^3, \;\;  \omega_3^4 = \omega_1^2.
$$
Thus by (40) relation (41) takes the form
\begin{equation}\label{eq:42}
d^2 A_0 \equiv (\lambda_{22} (\omega^2)^2 + 
2 \lambda_{23}  \omega^2 \omega^3 + \lambda_{33} 
(\omega^3)^2) A_4 +  2\omega^2  \omega^3  A_5 
     \pmod{T_x (V)}.
\end{equation}

Note that the coefficient in $A_5$ in equation (42) coincides 
with the first fundamental form (37) of 
a hypersurface $V \subset (M, c)$.

 Denote by $\widetilde{II}$ the coefficient in $A_4$ in 
equation (42):
$$
\widetilde{II} = \lambda_{22} (\omega^2)^2 
+ 2 \lambda_{23} \omega^2 \omega^3 + \lambda_{33} (\omega^3)^2. 
$$
Then equation (42) takes the form 
\begin{equation}\label{eq:43}
d^2 A_0 = \widetilde{II} A_4 +  I  A_5  \pmod{T_x (V)}.
\end{equation}

If we multiply expression (43) by a point 
$A_1 - x A_0$, then by (1) and (8), we find that 
$$
(d^2 A_0, A_1 - x A_0) = - (\widetilde{II} - x I).
$$
The expression in the parentheses of the right-hand side 
is a {\it pencil} of the second fundamental forms of 
a hypersurface $V \subset (M, c)$:

\begin{equation}\label{eq:44}
\widetilde{II} - x I = \lambda_{22} (\omega^2)^2 
+ 2 (\lambda_{23}  - x) \omega^2 \omega^3 
+ \lambda_{33} (\omega^3)^2.
\end{equation} 

The matrix of their coefficients is 
$$
({\widetilde{\lambda}}_{ab}) 
= \pmatrix{\lambda_{22} & \lambda_{23} -x\cr 
            \lambda_{23} -x  & \lambda_{33} \cr}. 
$$

From the pencil (43) we will take the form whose matrix is apolar 
to the matrix $(g_{ab})$, that is, the matrix satisfying the 
condition
\begin{equation}\label{eq:45}
{\widetilde{\lambda}}_{ab} g^{ab} = 0.
\end{equation}
Since by (38) we have
$$
 (g^{ab}) = \pmatrix{0 & 1 \cr
                     1 & 0 \cr},
$$
it follows that 
$$
{\widetilde{\lambda}}_{ab} g^{ab} =  2 (\lambda_{23} -x). 
$$
Thus condition (45) leads to the relation
\begin{equation}\label{eq:46}
x = \lambda_{23}.
\end{equation}

Condition (45) singles out from a pencil (44) of the second 
fundamental forms of a hypersurface $V \subset (M, c)$ 
a conformally invariant fundamental form 
\begin{equation}\label{eq:47}
II = \widetilde{II} - \lambda_{23} I 
= \lambda_{22} (\omega^2)^2  + \lambda_{33} (\omega^3)^2. 
\end{equation}
Its matrix has the form 
\begin{equation}\label{eq:48}
 (h_{ab}) = (\lambda_{ab}) - \lambda_{23} (g_{ab}) = 
\pmatrix{\lambda_{22} & 0 \cr
          0 & \lambda_{33} \cr}
\end{equation}
and is diagonal.

Consider singular points of the map $f (M^2 \times P^1) = V^3 \subset (M, c)$. We will look for these points in the form 
$X = A_1 - s A_0$. At these points the dimension 
of the tangent subspace $T_X (V)$ must be reduced. 
By (3), (11), and (34), we have 
\begin{equation}\label{eq:49}
dA_1 = \omega_1^0 A_0 + \omega_1^1 A_1 + \omega_1^2 A_2 
+ \omega_1^3 A_3. 
\end{equation}
Applying equations (49) and (35), we find that 
$$
d(A_1 - s A_0) = (\omega_1^0 - x \omega_0^0 - dx) A_0 
+ (\omega_1^1 - x \omega_0^1) A_1 + (\omega_1^2 - x \omega_0^2) 
A_2  + (\omega_1^3 - x \omega_0^3) A_3. 
$$
Further by (40) we obtain  
$$
d(A_1 - s A_0) \equiv ((\lambda_{23} - s) A_2 
+ \lambda_{22} A_3) \omega^2 + 
(\lambda_{33} A_2 + (\lambda_{23} - s) A_3 \omega^3)
(\pmod A_0, A_1).
$$
The tangent subspace $T_X (V)$ is determined by the points 
$A_0, A_1, (\lambda_{23} - s) A_2 + \lambda_{22} A_3$, and 
$\lambda_{33} A_2 + (\lambda_{23} - s) A_3$. Thus the dimension 
of the tangent subspace is reduced  only in the points  
$X = A_1 - s A_0$ in which
$$
\det{\pmatrix{\lambda_{23} - s & \lambda_{22} \cr
              \lambda_{33}  & \lambda_{23} - s \cr}} = 0.
$$
This equation can be written as 
\begin{equation}\label{eq:50}
s^2 - 2 \lambda_{23} s + (\lambda_{23}^2 - 
\lambda_{22} \lambda_{33}) = 0.  
\end{equation}
Denote by $s_1$ and $s_2$ the roots of 
this equation. They are calculated by the following formula: 
$$
s_{1, 2} = \lambda_{23} \pm \sqrt{\lambda_{22} \lambda_{33}}.
$$
The points $F_1 = A_1 - s_1 A_0$ and $F_2 = A_1 - s_2 A_0$ 
are singular points of an isotropic geodesic $l = A_0 A_1$ 
of a hypersurface $V$.

By Vieta's theorem, it follows from equation (50) 
that
$$
s_1 + s_2 = 2\lambda_{23}.
$$
Thus the point $H = A_1 - \lambda_{23} A_0$ is 
the fourth harmonic point $H$ of the point $A_0$ with respect 
to the points $F_1$ and $F_2$ on the line $l = A_0 A_1$. 
The singular points $F_1$ and $F_2$ are located symmetrically 
with respect to the points $A_0$ and $H$.

Now the conformally invariant second fundamental 
form $II$ of a hypersurface $V \subset (M, c)$ can be written as
$$
II = -  (d^2 A_0, H).
$$

We take a moving frame whose  vertex $A_1$ coincides with the point $H$.
This implies  $\lambda_{23} = 0$. As a result, equation (50) 
becomes 
$$
s^2 - h_{22} h_{33} = 0,  
$$  
and 
\begin{equation}\label{eq:51}
s_{1, 2} = \pm \sqrt{h_{22} h_{33}}.  
\end{equation}

The following theorem follows from relation (51).
\begin{theorem} 
\begin{description}
\item[(a)] The second 
fundamental form $II$ of a hypersurface $V \subset (M, c)$ 
at a point $A_0$ is positive definite or negative definite 
if and only if the isotropic geodesic $l = A_0 A_1$ 
through the point $x = A_0$ bears two real singular points. 
If at a point $x = A_0$ this form is an indeterminate form of 
rank two, then the singular points on the straight line 
$l = A_0 A_1$ are complex conjugate.

\item[(b)] The second 
fundamental form $II$ of a hypersurface $V \subset (M, c)$ 
at a point $x = A_0$ has the rank less than two if and only if
 the singular points on the isotropic geodesic $l = A_0 A_1$ 
coincide. In this case the point $H$ coincides with this multiple 
 singular point.
\end{description}
\end{theorem}

On a lightlike hypersurface $V$, 2-planes 
$A_0 \wedge  A_1 \wedge  A_2$ and $A_0 \wedge  A_1 \wedge  A_3$ 
of an isotropic 
frame bundle compose an $\alpha$- and $\beta$-distribution. 
Denote them by $\Delta_\alpha$ and $\Delta_\beta$.  These distributions are defined on $V$ by the equations
\begin{equation}\label{eq:52}
\omega^3 = 0 \hspace*{7mm}
 (\alpha)  \;\;\;\;\;\;\;\;\; \omega^2 = 0. \hspace*{7mm} 
(\beta)  
\end{equation}

In general, the distributions 
$\Delta_\alpha$ and $\Delta_\beta$  are 
not integrable. Let us find the conditions of their 
integrability. 

Exterior differentiation of equation (52$\alpha$) gives 
the following exterior quadratic equation 
$$
\omega^1 \wedge \omega_1^3 = 0.
$$
Substituting the value of the form $\omega_1^3$ from (40) into 
this equation and taking into account (48), we find 
that the distribution 
$\Delta_\alpha$ is integrable if and only if 
\begin{equation}\label{eq:53}
h_{22}  = 0.   
\end{equation}
Similarly the distribution 
$\Delta_\beta$ is integrable if and only if 
\begin{equation}\label{eq:54}
h_{33}  = 0.   
\end{equation}

Comparing the conditions (53) and (54) with relations (51) 
we arrive at the following result.

\begin{theorem} If at least one of the 
isotropic distributions 
$\Delta_\alpha$ and $\Delta_\beta$  on  a lightlike hypersurface 
$V \subset (M, c)$ is integrable, then the singular points on each of its isotropic generators coincide.
\end{theorem}

If both isotropic distributions 
$\Delta_\alpha$ and $\Delta_\beta$ are integrable on a 
hypersurface $V$, then conditions (53) and (54) are satisfied 
simultaneously, and  the second fundamental form $II$ of 
$V$ vanishes. But this means that hypersurface $V$ 
is totally umbilical. This implies the following result.

\begin{theorem} Both isotropic distributions 
$\Delta_\alpha$ and $\Delta_\beta$  on  a lightlike hypersurface 
$V \subset (M, c)$ are integrable if and only if the hypersurface 
$V$ is totally umbilical. 
\end{theorem}

\section{SOME PROPERTIES OF LIGHTLIKE \\ HYPERSURFACES}

We will pass now to the study of properties of 
a lightlike hypersurface $V \subset (M, c)$ connected with 
third- and higher-order differential neighborhoods.

We make a reduction in our isotropic second-order 
frame bundle by taking a specialized frame 
whose vertex  $A_1 \in l$ coincides 
with the fourth harmonic point $H$ of the point $A_0$ with 
respect to the singular points $F_1$ and $F_2$ of the straight line 
$l = A_0 A_1$. Then we obtain 
$$
\lambda_{23} = 0, \;\; h_{22} = \lambda_{22}, \;\; h_{33} = 
\lambda_{33}, 
$$ 
and equations (40) become 
\begin{equation}\label{eq:55}
\omega_1^3 = h_{22} \omega^2, \;\;
\omega_1^2 = h_{33} \omega^3.
\end{equation}

By (12), (14), (15), (18), and (19), exterior differentiation of 
equations (55) gives 
\begin{equation}\label{eq:56}
\left\{
\renewcommand{\arraystretch}{1.3}
\begin{array}{ll}
 \Delta h_{22}   \wedge  \omega^2 
+ (- \omega^0_1 + h_{22} h_{33} \omega^1 
- 2 a_1 \omega^2 - 2 b_1 \omega^3) \wedge \omega^3 = 0, \\ 
 (- \omega^0_1 + h_{22} h_{33} \omega^1 
- 2 a_1 \omega^2 - 2 b_1 \omega^3) \wedge \omega^2 
+ \Delta h_{33} \wedge  \omega^3 = 0, 
\end{array}
\renewcommand{\arraystretch}{1}
\right.
\end{equation}
where
\begin{equation}\label{eq:57}
\left\{
\renewcommand{\arraystretch}{1.3}
\begin{array}{ll}
\Delta h_{22} = dh_{22} + h_{22} (\omega^0_0 - 2 \omega^2_2 
- \omega^1_1) + 2 a_0 \omega^1, \\
\Delta h_{33} = d h_{33} + h_{33} (\omega^0_0 + 2 \omega^2_2 
- \omega^1_1) + 2 b_0 \omega^1.
\end{array}
\renewcommand{\arraystretch}{1}
\right.
\end{equation}
By Cartan's lemma, it follows from (56) that
\begin{equation}\label{eq:58}
\left\{
\renewcommand{\arraystretch}{1.3}
\begin{array}{ll}
\Delta h_{22}   = h_{222}  \omega^2 + h_{223} 
\omega^3,  \\ 
 \omega^0_1 =  h_{22} h_{33} \omega^1 -(h_{223} + 2a_1)  \omega^2  
- (h_{233} + 2 b_1)\omega^3, \\
\Delta h_{33}  = h_{233}  \omega^2  + h_{333} \omega^3.
\end{array}
\renewcommand{\arraystretch}{1}
\right.
\end{equation}

We will apply now equations (58) to totally umbilical 
hypersurfaces $V \subset (M, c)$. For such hypersurfaces we have $h_{22} = h_{33} = 0$. As a result,  equations (55) take the form
\begin{equation}\label{eq:59}
\omega_1^3 = 0, \;\;\omega_1^2 = 0,
\end{equation}
and equations (58) imply that 
\begin{equation}\label{eq:60}
a_0 = 0, \;\; b_0 = 0,
\end{equation}
\begin{equation}\label{eq:61}
h_{222} = h_{223} = h_{233} = h_{333}  = 0,
\end{equation}
and 
\begin{equation}\label{eq:62}
 \omega^0_1 =   - 2 (a_1  \omega^2   +  b_1\omega^3).
\end{equation}

Conditions (60) mean that the isotropic distributions 
$\Delta_\alpha$ and $\Delta_\beta$ are principal. 
Moreover, it follows now from equations (49) that 
\begin{equation}\label{eq:63}
d H = \omega^0_0 H - 2 (a_1  \omega^2   +  b_1\omega^3) A_0.
\end{equation}

This implies the following result.

\begin{theorem} A lightlike totally umbilical hypersurface 
$V \subset (M, c)$ possesses the following properties:
\begin{description}
\item[(a)] The isotropic distributions 
$\Delta_\alpha$ and $\Delta_\beta$ are integrable and principal.

 \item[(b)] The multiple singular point $H$ of the isotropic 
geodesic $l = A_0 A_1$ describes an isotropic line tangent to the 
straight line $l$ at the point $H$.

\item[(c)] If $a_1 = b_1 = 0$, then the point $H$ is fixed, 
and a totally umbilical hypersurface is an isotropic cone $C_H$ 
with vertex $H$.
\end{description}
\end{theorem}

{\sf Proof}. The statement (a) follows from the fact that 
on a hypersurface $V$, 
the isotropic distributions $\Delta_\alpha$ and $\Delta_\beta$ 
are defined by equations (52) and correspond to the values 
$\lambda = \infty$ and $\mu = \infty$ in equations (9) and (10). 
Hence for $a_0 = b_0 = 0$, these values satisfy equations (32)
 and (33) defining the principal isotropic  distributions. 

The statement (b)  follows immediately from equation (63).

Note  that the conditions $a_1 = b_1 = 0$ along with conditions 
(60) imply that   the values 
$\lambda = \infty$ and $\mu = \infty$ are multiple roots of 
equations (32) and (33). This implies that the statement 
(c) can be also formulated as follows:
\begin{description}
\item[(c')] A lightlike totally umbilical hypersurface 
$V \subset (M, c)$ is an isotropic cone if and only if it 
bears multiple isotropic 
distributions $\Delta_\alpha$ and $\Delta_\beta$.
\end{description}

Note also that in this case the integral surfaces of the  
distributions 
$\Delta_\alpha$ and $\Delta_\beta$ on a hypersurface $V$ 
are two-dimensional plane generators of the cone $C_H$.

\section{CONSTRUCTION OF A CANONICAL \\ 
DISTRIBUTION OF FRAMES FOR  \\ A
LIGHTLIKE HYPERSURFACE}

We associated a family of the second-order frames 
with a point $x = A_0$ of 
a lightlike totally umbilical hypersurface $V \subset (M, c)$ 
in such away that  the vertex $A_1$ coincides with the harmonic pole $H$ 
of the isotropic tangent $A_0 A_1$. But the points $A_2$ and 
$A_3$ of these frames can move freely in $\alpha$- 
and $\beta$-planes containing the straight line $A_0 A_4$, and 
its point $A_4$ can move freely along 
 the isotropic straight line $A_0 A_1$ that is conjugate 
to the screen subspace $A_0 \wedge A_2 \wedge A_3$ with respect to 
the isotropic cone $C_x$. 

For a fixed point $x = A_0$, by equations (16) and (40), 
we find that  
$$
\begin{array}{ll}
\renewcommand{\arraystretch}{1.3}
\delta A_2 = \pi_2^0 A_0 + \pi_2^1 A_1 + \pi_2^2 A_2,\\
\delta A_3 = \pi_3^0 A_0 + \pi_3^1 A_1 + \pi_3^3 A_3.
\renewcommand{\arraystretch}{1}
\end{array}
$$
Here $\pi_2^0, \pi_2^1, \pi_3^0$, and $\pi_3^1$ 
are fiber forms defining a displacement of the points 
$A_2$ and $A_3$ in the corresponding isotropic 2-planes. 

In order to find the points 
$A_2$ and $A_3$ uniquely in these 2-planes, one needs 
to make the above mentioned fiber forms vanish. 
However, this must be done in such a way that a 
fixing of the points 
$A_2$ and $A_3$ would be intrinsically connected with 
the geometry of a hypersurface $V$. The latter can be achieved 
by fixing in a certain way the coefficients $h_{abc}$  occurring 
in equations (58). These coefficients are associated with a third-order 
neighborhood of a hypersurface $V$. 

To this end, we take exterior derivatives of equations (58). 
As a result, we obtain the following exterior quadratic 
equations: 
\begin{equation}\label{eq:64}
\left\{
\renewcommand{\arraystretch}{1.3}
\begin{array}{ll}
\Delta h_{222} \wedge \omega^2 + \Delta h_{223} \wedge \omega^3 
+ H_{22} = 0, \\
\Delta h_{223} \wedge \omega^2 + \Delta h_{233} \wedge \omega^3 
+ H_{23} = 0, \\
\Delta h_{233} \wedge \omega^2 + \Delta h_{333} \wedge \omega^3 
+ H_{33} = 0, 
\end{array}
\renewcommand{\arraystretch}{1}
\right.
\end{equation}
where 
$$
\renewcommand{\arraystretch}{1.3}
\begin{array}{ll}
\Delta h_{222} = dh_{222} 
+ h_{222} (2\omega^0_0 - 3 \omega^2_2 - \omega^1_1) 
+ 2 a_0 \omega_2^1 + 3 h_{22} \omega_2^0 - 3 (h_{22})^2 
\omega_3^1,\\
\Delta h_{223} 
= dh_{223} + h_{223} (2\omega^0_0 - \omega^2_2 - \omega^1_1) 
+ 2 a_0 \omega_3^1 - h_{22} \omega_3^0 + h_{22} h_{33} 
\omega_2^1,  \\
\Delta h_{233} = 
dh_{233} + h_{233} (2\omega^0_0 + \omega^2_2 - \omega^1_1) 
+ 2 b_0 \omega_2^1 - h_{33} \omega_2^0 +  h_{22} h_{33} 
\omega_3^1,\\
\Delta h_{333} 
= dh_{333} + h_{333} (2\omega^0_0 + 3 \omega^2_2 - \omega^1_1) 
+ 2 b_0 \omega_3^1  + 3 h_{33} \omega_3^0 
- 3 (h_{33})^2 \omega_2^1.
\end{array}
\renewcommand{\arraystretch}{1}
$$
and $H_{22}, H_{23}$, and $H_{33}$ are 2-forms that 
are linear combinations of the products 
$\omega^1 \wedge \omega^2, \; \omega^2 \wedge \omega^3$, 
and $\omega^1 \wedge \omega^3$ of the basis forms 
$\omega^1, \omega^2$, and $\omega^3$. 

Equations (64) imply that the 1-forms $\Delta h_{222}, \;
\Delta h_{223},\; \Delta h_{233}$, and $\Delta h_{333}$ are 
linear combinations of  the basis forms 
$\omega^1, \omega^2$, and $\omega^3$. 

For a fixed point $x = A_0$, i.e., for $\omega^1 = \omega^2 = 
 \omega^3 = 0$, these forms vanish, and their expressions become
\begin{equation}\label{eq:65}
\renewcommand{\arraystretch}{1.3}
\begin{array}{ll}
\Delta_\delta h_{222} =  \delta h_{222} 
+ h_{222} (2\pi^0_0 - 3 \pi^2_2 - \pi^1_1) 
+ 2 a_0 \pi_2^1 + 3 h_{22} \pi_2^0 - 3 (h_{22})^2 
\pi_3^1 = 0,\\
\Delta_\delta h_{223} 
= \delta h_{223} + h_{223} (2\pi^0_0 - \pi^2_2 - \pi^1_1) 
+ 2 a_0 \pi_3^1 - h_{22} \pi_3^0 + h_{22} h_{33} 
\pi_2^1 = 0,  \\
\Delta_\delta h_{233} = 
\delta h_{233} + h_{233} (2\pi^0_0 + \pi^2_2 - \pi^1_1) 
+ 2 b_0 \pi_2^1 - h_{33} \pi_2^0 +  h_{22} h_{33} 
\pi_3^1 = 0,\\
\Delta_\delta h_{333} 
= \delta h_{333} + h_{333} (2\pi^0_0 + 3 \pi^2_2 - \pi^1_1) 
 + 2 b_0 \pi_3^1  + 3 h_{33} \pi_3^0 - 3 (h_{33})^2 \pi_2^1 = 0.
\end{array}
\renewcommand{\arraystretch}{1}
\end{equation}

Equations (65) contain the fiber forms $\pi_2^0, \; \pi_2^1, 
\pi_3^0$, and $\pi_3^1$ defining a displacement of the points 
$A_2$ and $A_3$ in the $\alpha$- and $\beta$-planes 
$A_0 \wedge A_1 \wedge A_2$ and $A_0 \wedge A_1 \wedge A_3$. 
Consider the determinant $D$ of the matrix 
of coefficients in these fiber forms in equations (65):
$$
D = \det \pmatrix{3 h_{22} & 2 a_0 & 0 & - 3 (h_{22})^2 \cr 
              0        &h_{22} h_{33} & - h_{22} & 2 a_0 \cr
              -  h_{33} & 2 b_0 & 0 & 3 h_{22} h_{33} \cr
              0        &- 3 (h_{33})^2 & 3 h_{33} & 2 b_0 \cr}.
$$
Calculating this determinant, we find that 
\begin{equation}\label{eq:66}
D = 4 (3 h_{22} b_0 + h_{33} a_0)(h_{22} b_0 + 3 h_{33} a_0).
\end{equation}

If this determinant does not vanish, $D \neq 0$, then 
equations (65) imply that the quantities $h_{222}, \; 
h_{223}, \; h_{233}$, and $h_{333},$ occurring in equations (58) can be 
simultaneously reduced to 0 by means of the 
 fiber forms $\pi_2^0, \pi_2^1, \pi_3^0$, and $\pi_3^1$ 
(see [O 62]). As a result, 
the points $A_2$ and $A_3$ are uniquely  determined 
in the planes  $\alpha = A_0 \wedge A_1 \wedge A_2$ and 
$\beta = A_0 \wedge A_1 \wedge A_3$, and we arrive at 
a family of  third-order moving frames associated 
with a point $x = A_0 \in V \subset (M, c)$.

 With respect to a third-order frame 
we have constructed, equations (58) take the form 
\begin{equation}\label{eq:67}
\left\{
\renewcommand{\arraystretch}{1.3}
\begin{array}{ll}
d h_{22}   + h_{22} (\omega_0^0 - 2 \omega_2^2 
- \omega_1^1)   = - 2 a_0   \omega^1,  \\ 
 \omega^0_1 =  h_{22} h_{33} \omega^1 - 2a_1  \omega^2  
- 2 b_1 \omega^3, \\
dh_{33} + h_{33} (\omega_0^0 + 2 \omega_2^2 
- \omega_1^1)   = - 2 b_0   \omega^1.
\end{array}
\renewcommand{\arraystretch}{1}
\right.
\end{equation}
This proves the following result. 

\begin{theorem} If on a lightlike hypersurface $V$ 
the determinant $D$ does not vanish, then it is possible to 
construct a third-order frame bundle on $V$ intrinsically 
connected with the geometry of $V$. In this frame bundle, 
$h_{abc} = 0$. 
\end{theorem}

Note that if the $CO (2, 2)$-structure on 
a manifold $(M, c)$ is conformally flat, then a third-order 
frame bundle indicated above cannot be constructed since 
for a conformally flat structure we have $a_0 = b_0 = 0$, 
and consequently, $D = 0$. However, for 
a conformally semiflat $CO (2, 2)$-structure the above 
construction is possible. A construction of 
a canonical frame bundle for lightlike totally 
umbilical hypersurfaces is also impossible since for 
them $h_{22} = h_{33} = 0$, and consequently, $D = 0$.

In order to complete our construction of a canonical 
frame bundle, we also have to  fix the vertex $A_4$ on 
the isotropic straight line $A_0 A_4$ which is conjugate to 
the screen subspace $S_x = A_0 \wedge A_2 \wedge A_3$ 
with respect to the isotropic cone $C_x$. 
 This can be 
done in the same way as we did 
for a lightlike hypersurface $V \subset (M, c)$ whose conformal 
structure $c$ is of Lorentzian signature, $c = CO (n-1, 1)$.  
The family of straight lines $A_0 A_4$ associated with 
a hypersurface $V$ is an isotropic congruence (see [AG 99b]) 
each ray of which bears two singular points $F_1^\prime$ and 
$F_2^\prime$. To complete our specialization of 
moving frames, we choose a frame whose vertex $A_4$ coincides with the 
harmonic pole $H'$ of the point $A_0$ with respect to 
singular points $F_1^\prime$ and $F_2^\prime$ (see Figure 3). 
Since the singular points are defined in a fourth-order 
differential 
neighborhood of a point $x \in V$, the point $A_4$ is defined 
also in this neighborhood.

Thus we arrive at the following result.

\begin{theorem} If $D \neq 0$, a canonical frame bundle 
on a lightlike hypersurface $V \subset (M, c), c = CO (2, 2)$ 
is defined by elements of a fourth-order differential 
neighborhood of a point $x \in V$.
\end{theorem}

\vspace*{2.8in}
\begin{center}
Figure 3
\end{center}

\noindent
{\em Authors' addresses}:\\

\noindent
\begin{tabular}{ll}
M. A. Akivis &                           V. V. Goldberg \\
Department of Mathematics      &      Department of Mathematics\\
Jerusalem College of Technology  & 
                       New Jersey Institute of Technology \\
-- Mahon Lev, P. O. B. 16031 & University Heights \\
Jerusalem 91160, Israel & Newark, NJ 07102, U.S.A. \\
\\
E-mail address: akivis@avoda.jct.ac.il & E-mail address: 
                                        vlgold@numerics.njit.edu    
\end{tabular}

\end{document}